\documentclass[11pt]{article}
\usepackage{amsmath}
\usepackage{amsfonts}
\usepackage{amssymb}
\usepackage{amsthm}
\newtheorem*{lem}{Lemma 6.3}
\begin{document}

\title{Proof of Lemma 6.3 in ``The crossing number of $K_{4,n}$ on the
torus and the Klein bottle"}
\author{Pak Tung Ho}

\maketitle

The Lemma 6.3 in ``The crossing number of $K_{4,n}$ on the torus and
the Klein bottle" is false, which says that

\begin{lem}\label{thankgod}
For any drawing $D$ of $K_{4,n}$ on the Klein bottle, let $A$ be the
matrix defined by $A_{ij}=\widetilde{cr_D}(a_i,a_j)$. Then it is
impossible for the following to hold for some distinct $i_j$, $1\leq
j\leq 5$:
\begin{eqnarray}\label{forbid} \left(
\begin{array}{ccccc}
    {A_{i_1i_1}} & {A_{i_1i_2}}& {A_{i_1i_3}} & {A_{i_1i_4}}& {A_{i_1i_5}}\\
    {A_{i_2i_1}} & {A_{i_2i_2}}& {A_{i_2i_3}} & {A_{i_2i_4}}& {A_{i_2i_5}}\\
    {A_{i_3i_1}} & {A_{i_3i_2}}& {A_{i_3i_3}} & {A_{i_3i_4}}& {A_{i_3i_5}}\\
    {A_{i_4i_1}} & {A_{i_4i_2}}& {A_{i_4i_3}} & {A_{i_4i_4}}&
    {A_{i_4i_5}}\\
    {A_{i_5i_1}} & {A_{i_5i_2}}& {A_{i_5i_3}} & {A_{i_5i_4}}&
    {A_{i_5i_5}}
    \end{array} \right)
= \left(
\begin{array}{ccccc}
    {0} & {0} & {0} & {0}& {0}\\
    {0} & {0} & {0} & {0}& {0}\\
    {0} & {0} & {0}& {1}& {1}\\
    {0} & {0} & {1}& {0}& {1}\\
    {0} & {0} & {1}& {1}& {0}
    \end{array} \right).\end{eqnarray}

\end{lem}

For the counterexample, one can refer to the following drawing of
$K_{4,5}$ on the Klein bottle with $a_j=a_{i_j}$, where
$\blacktriangle$ represents the vertices $a_3$, $a_4$ and $a_5$.

\begin{center}
\begin{picture}(0,0)(50,104)

\put(-2,8){$\blacksquare$}

\put(-2,98){$\blacksquare$}

\put(98,8){$\blacksquare$}

\put(98,98){$\blacksquare$}

\put(-2,-2){$b_1$}

\put(-2,108){$b_1$}

\put(98,-2){$b_1$}

\put(98,108){$b_1$}

\put(23,7){$\bullet$}  \put(18,0){$a_1$}

\put(48,8){$\blacksquare$} \put(43,3){$b_2$}

\put(73,7){$\bullet$} \put(81,0){$a_2$}

\put(23,97){$\bullet$} \put(15,102){$a_1$}

\put(48,98){$\blacksquare$} \put(57,103){$b_2$}

\put(73,97){$\bullet$} \put(82,102){$a_2$}

\put(0,10){\line(1,0){100}} \put(0,100){\line(1,0){100}}

\multiput(0,10)(0,4){22}{$\cdot$}
\multiput(100,10)(0,4){22}{$\cdot$}

\put(-2,45){$\wedge$}
 \put(98,45){$\vee$}

\put(25,100){\line(1,-1){25}} \put(75,100){\line(-1,-1){25}}
\put(47,73){$\blacksquare$} \put(47,83){$b_4$}

\put(25,10){\line(-1,1){25}}

\put(-4,33){$\blacksquare$} \put(-14,30){$b_3$}

 \put(97,53){$\blacksquare$}  \put(107,57){$b_3$}
\put(77,100){\line(1,-2){21}}

\put(47,43){$\blacktriangle$} \put(51,45){\line(0,1){29}}
\put(51,45){\line(0,-1){32}} \put(51,45){\line(5,1){52}}
\put(51,45){\line(3,-2){52}}

\put(75,40){$\blacktriangle$} \put(79,42){\line(-3,4){27}}
\put(79,41){\line(3,2){20}} \put(79,42){\line(3,-4){23}}
\put(79,42){\line(3,-1){23}} \put(0,60){\line(1,-1){51}}

\put(38,33){$\blacktriangle$} \put(42,35){\line(1,-3){8}}
\put(42,35){\line(5,-2){60}} \put(42,35){\line(1,5){8}}
\put(42,35){\line(-1,0){40}} \put(42,37){\line(-2,3){42}}

 \put(30,-18){Figure 1.}

\end{picture}
\end{center}
$$$$\\\\\\\\\\\\

\begin{center}
APPENDIX. The possible drawings of the $K_{2,4}$ containing $a_1$
and $a_2$ in torus which has no crossings.
\end{center}

To obtain these drawings, we combine all the possibilities of
different ways of drawing edges $a_1b_3$, $a_2b_3$, $a_1b_4$ and
$a_2b_4$, as in Figure 6(i) to 6(viii) in \cite{Tung}. Figure
2(s)(t) was obtained by drawing the pair $a_1b_3$, $a_2b_3$ as in
Figure 6(s) in \cite{Tung} and the pair $a_1b_4$, $a_2b_4$ as in
Figure 6(t) in \cite{Tung} where s, t equals i, ii,...., viii. For
some s, t, Figure 2(s)(t) is omitted since it has crossings.
Moreover, since Figure 2(s)(t) and Figure 2(t)(s) are identical,
only Figure 2(s)(t), where s$\leq$t, is shown.\\\\

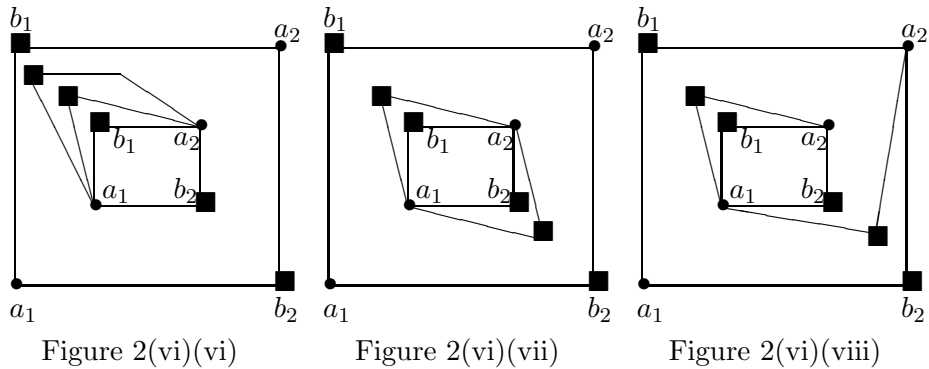
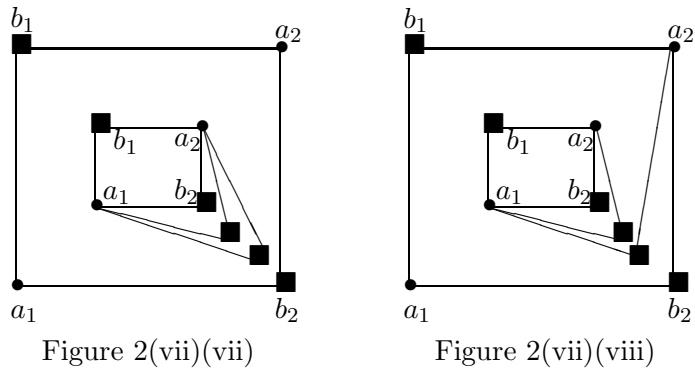
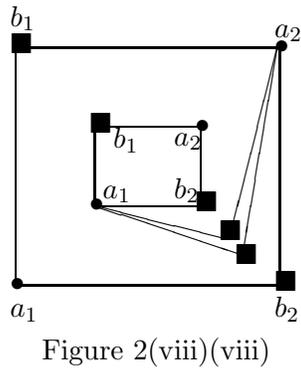
\begin{center}
\begin{picture}(0,0)(160,104)

\put(-2,8){$\bullet$} \put(-2,-2){$a_1$}

\put(-2,98){$\blacksquare$} \put(-2,108){$b_1$}

\put(98,8){$\blacksquare$} \put(98,-2){$b_2$}

\put(98,98){$\bullet$} \put(98,104){$a_2$}

\put(0,10){\line(1,0){100}} \put(0,10){\line(0,1){90}}
\put(0,100){\line(1,0){100}} \put(100,10){\line(0,1){90}}

\put(28,38){$\bullet$} \put(33,43){$a_1$}

\put(68,38){$\blacksquare$} \put(60,43){$b_2$}

\put(28,68){$\blacksquare$} \put(37,62){$b_1$}

\put(68,68){$\bullet$} \put(60,63){$a_2$}

 \put(30,40){\line(1,0){40}}  \put(30,40){\line(0,1){30}}
 \put(30,70){\line(1,0){40}} \put(70,40){\line(0,1){30}}

  \put(30,40){\line(-1,2){21}} \put(100,100){\line(-6,-1){90}}
 \put(7,80){$\blacksquare$}

 \put(100,100){\line(-4,-1){80}} \put(30,40){\line(-1,4){10}}
 \put(17,75){$\blacksquare$}

 \put(20,-18){Figure 2(i)(i)}

\end{picture}
\begin{picture}(0,0)(50,104)

\put(-2,8){$\bullet$} \put(-2,-2){$a_1$}

\put(-2,98){$\blacksquare$} \put(-2,108){$b_1$}

\put(98,8){$\blacksquare$} \put(98,-2){$b_2$}

\put(98,98){$\bullet$} \put(98,104){$a_2$}

\put(0,10){\line(1,0){100}} \put(0,10){\line(0,1){90}}
\put(0,100){\line(1,0){100}} \put(100,10){\line(0,1){90}}

\put(28,38){$\bullet$} \put(33,43){$a_1$}

\put(68,38){$\blacksquare$} \put(60,43){$b_2$}

\put(28,68){$\blacksquare$} \put(37,62){$b_1$}

\put(68,68){$\bullet$} \put(60,63){$a_2$}

 \put(30,40){\line(1,0){40}}  \put(30,40){\line(0,1){30}}
 \put(30,70){\line(1,0){40}} \put(70,40){\line(0,1){30}}

 \put(0,10){\line(6,1){90}}
 \put(70,70){\line(1,-2){20}} \put(85,25){$\blacksquare$}

  \put(30,40){\line(-1,2){21}} \put(100,100){\line(-6,-1){90}}
 \put(7,80){$\blacksquare$}

 \put(20,-18){Figure 2(i)(ii)}

\end{picture}
\begin{picture}(0,0)(-65,104)

\put(-2,8){$\bullet$} \put(-2,-2){$a_1$}

\put(-2,98){$\blacksquare$} \put(-2,108){$b_1$}

\put(98,8){$\blacksquare$} \put(98,-2){$b_2$}

\put(98,98){$\bullet$} \put(98,104){$a_2$}

\put(0,10){\line(1,0){100}} \put(0,10){\line(0,1){90}}
\put(0,100){\line(1,0){100}} \put(100,10){\line(0,1){90}}

\put(28,38){$\bullet$} \put(33,43){$a_1$}

\put(68,38){$\blacksquare$} \put(60,43){$b_2$}

\put(28,68){$\blacksquare$} \put(37,62){$b_1$}

\put(68,68){$\bullet$} \put(60,63){$a_2$}

 \put(30,40){\line(1,0){40}}  \put(30,40){\line(0,1){30}}
 \put(30,70){\line(1,0){40}} \put(70,40){\line(0,1){30}}

  \put(30,40){\line(-1,2){21}} \put(100,100){\line(-6,-1){90}}
 \put(7,80){$\blacksquare$}

 \put(0,10){\line(6,1){90}}
 \put(100,100){\line(-1,-6){12}} \put(85,25){$\blacksquare$}

 \put(20,-18){Figure 2(i)(iv)}

\end{picture}
\end{center}
$$$$\\\\\\\\\\\\

\begin{center}
\begin{picture}(0,0)(160,104)

\put(-2,8){$\bullet$} \put(-2,-2){$a_1$}

\put(-2,98){$\blacksquare$} \put(-2,108){$b_1$}

\put(98,8){$\blacksquare$} \put(98,-2){$b_2$}

\put(98,98){$\bullet$} \put(98,104){$a_2$}

\put(0,10){\line(1,0){100}} \put(0,10){\line(0,1){90}}
\put(0,100){\line(1,0){100}} \put(100,10){\line(0,1){90}}

\put(28,38){$\bullet$} \put(33,43){$a_1$}

\put(68,38){$\blacksquare$} \put(60,43){$b_2$}

\put(28,68){$\blacksquare$} \put(37,62){$b_1$}

\put(68,68){$\bullet$} \put(60,63){$a_2$}

 \put(30,40){\line(1,0){40}}  \put(30,40){\line(0,1){30}}
 \put(30,70){\line(1,0){40}} \put(70,40){\line(0,1){30}}

 \put(100,100){\line(-4,-1){80}} \put(30,40){\line(-1,4){10}}
 \put(17,75){$\blacksquare$}

  \put(0,10){\line(1,6){12}} \put(100,100){\line(-6,-1){90}}
 \put(7,80){$\blacksquare$}

 \put(20,-18){Figure 2(i)(v)}

\end{picture}
\begin{picture}(0,0)(50,104)

\put(-2,8){$\bullet$} \put(-2,-2){$a_1$}

\put(-2,98){$\blacksquare$} \put(-2,108){$b_1$}

\put(98,8){$\blacksquare$} \put(98,-2){$b_2$}

\put(98,98){$\bullet$} \put(98,104){$a_2$}

\put(0,10){\line(1,0){100}} \put(0,10){\line(0,1){90}}
\put(0,100){\line(1,0){100}} \put(100,10){\line(0,1){90}}

\put(28,38){$\bullet$} \put(33,43){$a_1$}

\put(68,38){$\blacksquare$} \put(60,43){$b_2$}

\put(28,68){$\blacksquare$} \put(37,62){$b_1$}

\put(68,68){$\bullet$} \put(60,63){$a_2$}

 \put(30,40){\line(1,0){40}}  \put(30,40){\line(0,1){30}}
 \put(30,70){\line(1,0){40}} \put(70,40){\line(0,1){30}}

  \put(30,40){\line(-1,2){21}} \put(100,100){\line(-6,-1){90}}
 \put(7,80){$\blacksquare$}

\put(30,40){\line(-1,4){10}} \put(70,70){\line(-4,1){50}}
\put(16,78){$\blacksquare$}

 \put(20,-18){Figure 2(i)(vi)}

\end{picture}
\begin{picture}(0,0)(-65,104)

\put(-2,8){$\bullet$} \put(-2,-2){$a_1$}

\put(-2,98){$\blacksquare$} \put(-2,108){$b_1$}

\put(98,8){$\blacksquare$} \put(98,-2){$b_2$}

\put(98,98){$\bullet$} \put(98,104){$a_2$}

\put(0,10){\line(1,0){100}} \put(0,10){\line(0,1){90}}
\put(0,100){\line(1,0){100}} \put(100,10){\line(0,1){90}}

\put(28,38){$\bullet$} \put(33,43){$a_1$}

\put(68,38){$\blacksquare$} \put(60,43){$b_2$}

\put(28,68){$\blacksquare$} \put(37,62){$b_1$}

\put(68,68){$\bullet$} \put(60,63){$a_2$}

 \put(30,40){\line(1,0){40}}  \put(30,40){\line(0,1){30}}
 \put(30,70){\line(1,0){40}} \put(70,40){\line(0,1){30}}

  \put(30,40){\line(-1,2){21}} \put(100,100){\line(-6,-1){90}}
 \put(7,80){$\blacksquare$}

\put(30,40){\line(4,-1){50}} \put(71,70){\line(1,-4){10}}
\put(77,27){$\blacksquare$}

 \put(20,-18){Figure 2(i)(vii)}

\end{picture}
\end{center}
$$$$\\\\\\\\\\\\

\newpage

\begin{center}
\begin{picture}(0,0)(50,104)

\put(-2,8){$\bullet$} \put(-2,-2){$a_1$}

\put(-2,98){$\blacksquare$} \put(-2,108){$b_1$}

\put(98,8){$\blacksquare$} \put(98,-2){$b_2$}

\put(98,98){$\bullet$} \put(98,104){$a_2$}

\put(0,10){\line(1,0){100}} \put(0,10){\line(0,1){90}}
\put(0,100){\line(1,0){100}} \put(100,10){\line(0,1){90}}

\put(28,38){$\bullet$} \put(33,43){$a_1$}

\put(68,38){$\blacksquare$} \put(60,43){$b_2$}

\put(28,68){$\blacksquare$} \put(37,62){$b_1$}

\put(68,68){$\bullet$} \put(60,63){$a_2$}

 \put(30,40){\line(1,0){40}}  \put(30,40){\line(0,1){30}}
 \put(30,70){\line(1,0){40}} \put(70,40){\line(0,1){30}}

  \put(30,40){\line(-1,2){21}} \put(100,100){\line(-6,-1){90}}
 \put(7,80){$\blacksquare$}

 \put(100,100){\line(-1,-6){12}} \put(85,25){$\blacksquare$}
 \put(87,30){\line(-6,1){57}}

 \put(20,-18){Figure 2(i)(viii)}

\end{picture}
\end{center}
$$$$\\\\\\\\\\\\

\begin{center}
\begin{picture}(0,0)(160,104)

\put(-2,8){$\bullet$} \put(-2,-2){$a_1$}

\put(-2,98){$\blacksquare$} \put(-2,108){$b_1$}

\put(98,8){$\blacksquare$} \put(98,-2){$b_2$}

\put(98,98){$\bullet$} \put(98,104){$a_2$}

\put(0,10){\line(1,0){100}} \put(0,10){\line(0,1){90}}
\put(0,100){\line(1,0){100}} \put(100,10){\line(0,1){90}}

\put(28,38){$\bullet$} \put(33,43){$a_1$}

\put(68,38){$\blacksquare$} \put(60,43){$b_2$}

\put(28,68){$\blacksquare$} \put(37,62){$b_1$}

\put(68,68){$\bullet$} \put(60,63){$a_2$}

 \put(30,40){\line(1,0){40}}  \put(30,40){\line(0,1){30}}
 \put(30,70){\line(1,0){40}} \put(70,40){\line(0,1){30}}

 \put(0,10){\line(6,1){90}}
 \put(70,70){\line(1,-2){20}} \put(85,25){$\blacksquare$}

 \put(0,10){\line(4,1){80}} \put(70,70){\line(1,-4){10}} \put(76,28){$\blacksquare$}

 \put(20,-18){Figure 2(ii)(ii)}

\end{picture}
\begin{picture}(0,0)(50,104)

\put(-2,8){$\bullet$} \put(-2,-2){$a_1$}

\put(-2,98){$\blacksquare$} \put(-2,108){$b_1$}

\put(98,8){$\blacksquare$} \put(98,-2){$b_2$}

\put(98,98){$\bullet$} \put(98,104){$a_2$}

\put(0,10){\line(1,0){100}} \put(0,10){\line(0,1){90}}
\put(0,100){\line(1,0){100}} \put(100,10){\line(0,1){90}}

\put(28,38){$\bullet$} \put(33,43){$a_1$}

\put(68,38){$\blacksquare$} \put(60,43){$b_2$}

\put(28,68){$\blacksquare$} \put(37,62){$b_1$}

\put(68,68){$\bullet$} \put(60,63){$a_2$}

 \put(30,40){\line(1,0){40}}  \put(30,40){\line(0,1){30}}
 \put(30,70){\line(1,0){40}} \put(70,40){\line(0,1){30}}

  \put(0,10){\line(1,6){12}} \put(70,70){\line(-4,1){60}}
 \put(7,80){$\blacksquare$}

 \put(0,10){\line(6,1){90}}
 \put(70,70){\line(1,-2){20}} \put(85,25){$\blacksquare$}

 \put(10,-18){Figure 2(ii)(iii)}

\end{picture}
\begin{picture}(0,0)(-65,104)

\put(-2,8){$\bullet$} \put(-2,-2){$a_1$}

\put(-2,98){$\blacksquare$} \put(-2,108){$b_1$}

\put(98,8){$\blacksquare$} \put(98,-2){$b_2$}

\put(98,98){$\bullet$} \put(98,104){$a_2$}

\put(0,10){\line(1,0){100}} \put(0,10){\line(0,1){90}}
\put(0,100){\line(1,0){100}} \put(100,10){\line(0,1){90}}

\put(28,38){$\bullet$} \put(33,43){$a_1$}

\put(68,38){$\blacksquare$} \put(60,43){$b_2$}

\put(28,68){$\blacksquare$} \put(37,62){$b_1$}

\put(68,68){$\bullet$} \put(60,63){$a_2$}

 \put(30,40){\line(1,0){40}}  \put(30,40){\line(0,1){30}}
 \put(30,70){\line(1,0){40}} \put(70,40){\line(0,1){30}}

 \put(0,10){\line(4,1){80}} \put(70,70){\line(1,-4){10}} \put(76,28){$\blacksquare$}

 \put(0,10){\line(6,1){90}}
 \put(100,100){\line(-1,-6){12}} \put(85,25){$\blacksquare$}

 \put(20,-18){Figure 2(ii)(iv)}

\end{picture}
\end{center}
$$$$\\\\\\\\\\\\

\begin{center}
\begin{picture}(0,0)(165,104)

\put(-2,8){$\bullet$} \put(-2,-2){$a_1$}

\put(-2,98){$\blacksquare$} \put(-2,108){$b_1$}

\put(98,8){$\blacksquare$} \put(98,-2){$b_2$}

\put(98,98){$\bullet$} \put(98,104){$a_2$}

\put(0,10){\line(1,0){100}} \put(0,10){\line(0,1){90}}
\put(0,100){\line(1,0){100}} \put(100,10){\line(0,1){90}}

\put(28,38){$\bullet$} \put(33,43){$a_1$}

\put(68,38){$\blacksquare$} \put(60,43){$b_2$}

\put(28,68){$\blacksquare$} \put(37,62){$b_1$}

\put(68,68){$\bullet$} \put(60,63){$a_2$}

 \put(30,40){\line(1,0){40}}  \put(30,40){\line(0,1){30}}
 \put(30,70){\line(1,0){40}} \put(70,40){\line(0,1){30}}

 \put(0,10){\line(4,1){80}} \put(70,70){\line(1,-4){10}} \put(76,28){$\blacksquare$}

  \put(0,10){\line(1,6){12}} \put(100,100){\line(-6,-1){90}}
 \put(7,80){$\blacksquare$}

 \put(20,-18){Figure 2(ii)(v)}

\end{picture}
\begin{picture}(0,0)(50,104)

\put(-2,8){$\bullet$} \put(-2,-2){$a_1$}

\put(-2,98){$\blacksquare$} \put(-2,108){$b_1$}

\put(98,8){$\blacksquare$} \put(98,-2){$b_2$}

\put(98,98){$\bullet$} \put(98,104){$a_2$}

\put(0,10){\line(1,0){100}} \put(0,10){\line(0,1){90}}
\put(0,100){\line(1,0){100}} \put(100,10){\line(0,1){90}}

\put(28,38){$\bullet$} \put(33,43){$a_1$}

\put(68,38){$\blacksquare$} \put(60,43){$b_2$}

\put(28,68){$\blacksquare$} \put(37,62){$b_1$}

\put(68,68){$\bullet$} \put(60,63){$a_2$}

 \put(30,40){\line(1,0){40}}  \put(30,40){\line(0,1){30}}
 \put(30,70){\line(1,0){40}} \put(70,40){\line(0,1){30}}

 \put(0,10){\line(4,1){80}} \put(70,70){\line(1,-4){10}} \put(76,28){$\blacksquare$}

\put(30,40){\line(-1,4){10}} \put(70,70){\line(-4,1){50}}
\put(16,78){$\blacksquare$}

 \put(20,-18){Figure 2(ii)(vi)}

\end{picture}
\begin{picture}(0,0)(-65,104)

\put(-2,8){$\bullet$} \put(-2,-2){$a_1$}

\put(-2,98){$\blacksquare$} \put(-2,108){$b_1$}

\put(98,8){$\blacksquare$} \put(98,-2){$b_2$}

\put(98,98){$\bullet$} \put(98,104){$a_2$}

\put(0,10){\line(1,0){100}} \put(0,10){\line(0,1){90}}
\put(0,100){\line(1,0){100}} \put(100,10){\line(0,1){90}}

\put(28,38){$\bullet$} \put(33,43){$a_1$}

\put(68,38){$\blacksquare$} \put(60,43){$b_2$}

\put(28,68){$\blacksquare$} \put(37,62){$b_1$}

\put(68,68){$\bullet$} \put(60,63){$a_2$}

 \put(30,40){\line(1,0){40}}  \put(30,40){\line(0,1){30}}
 \put(30,70){\line(1,0){40}} \put(70,40){\line(0,1){30}}

 \put(0,10){\line(6,1){90}}
 \put(70,70){\line(1,-2){20}} \put(85,25){$\blacksquare$}

\put(30,40){\line(4,-1){50}} \put(71,70){\line(1,-4){10}}
\put(77,27){$\blacksquare$}

 \put(20,-18){Figure 2(ii)(vii)}

\end{picture}
\end{center}
$$$$\\\\\\\\\\\\

\newpage

\begin{center}
\begin{picture}(0,0)(165,104)

\put(-2,8){$\bullet$} \put(-2,-2){$a_1$}

\put(-2,98){$\blacksquare$} \put(-2,108){$b_1$}

\put(98,8){$\blacksquare$} \put(98,-2){$b_2$}

\put(98,98){$\bullet$} \put(98,104){$a_2$}

\put(0,10){\line(1,0){100}} \put(0,10){\line(0,1){90}}
\put(0,100){\line(1,0){100}} \put(100,10){\line(0,1){90}}

\put(28,38){$\bullet$} \put(33,43){$a_1$}

\put(68,38){$\blacksquare$} \put(60,43){$b_2$}

\put(28,68){$\blacksquare$} \put(37,62){$b_1$}

\put(68,68){$\bullet$} \put(60,63){$a_2$}

 \put(30,40){\line(1,0){40}}  \put(30,40){\line(0,1){30}}
 \put(30,70){\line(1,0){40}} \put(70,40){\line(0,1){30}}

  \put(0,10){\line(1,6){13}} \put(70,70){\line(-3,1){60}}
 \put(8,86){$\blacksquare$}

\put(70,70){\line(-5,1){50}} \put(17,76){$\blacksquare$}
\put(0,10){\line(1,4){18}}

 \put(20,-18){Figure 2(iii)(iii)}

\end{picture}
\begin{picture}(0,0)(50,104)

\put(-2,8){$\bullet$} \put(-2,-2){$a_1$}

\put(-2,98){$\blacksquare$} \put(-2,108){$b_1$}

\put(98,8){$\blacksquare$} \put(98,-2){$b_2$}

\put(98,98){$\bullet$} \put(98,104){$a_2$}

\put(0,10){\line(1,0){100}} \put(0,10){\line(0,1){90}}
\put(0,100){\line(1,0){100}} \put(100,10){\line(0,1){90}}

\put(28,38){$\bullet$} \put(33,43){$a_1$}

\put(68,38){$\blacksquare$} \put(60,43){$b_2$}

\put(28,68){$\blacksquare$} \put(37,62){$b_1$}

\put(68,68){$\bullet$} \put(60,63){$a_2$}

 \put(30,40){\line(1,0){40}}  \put(30,40){\line(0,1){30}}
 \put(30,70){\line(1,0){40}} \put(70,40){\line(0,1){30}}

  \put(0,10){\line(1,6){13}} \put(70,70){\line(-3,1){60}}
 \put(8,86){$\blacksquare$}

 \put(0,10){\line(6,1){90}}
 \put(100,100){\line(-1,-6){12}} \put(85,25){$\blacksquare$}

 \put(20,-18){Figure 2(iii)(iv)}

\end{picture}
\begin{picture}(0,0)(-65,104)

\put(-2,8){$\bullet$} \put(-2,-2){$a_1$}

\put(-2,98){$\blacksquare$} \put(-2,108){$b_1$}

\put(98,8){$\blacksquare$} \put(98,-2){$b_2$}

\put(98,98){$\bullet$} \put(98,104){$a_2$}

\put(0,10){\line(1,0){100}} \put(0,10){\line(0,1){90}}
\put(0,100){\line(1,0){100}} \put(100,10){\line(0,1){90}}

\put(28,38){$\bullet$} \put(33,43){$a_1$}

\put(68,38){$\blacksquare$} \put(60,43){$b_2$}

\put(28,68){$\blacksquare$} \put(37,62){$b_1$}

\put(68,68){$\bullet$} \put(60,63){$a_2$}

 \put(30,40){\line(1,0){40}}  \put(30,40){\line(0,1){30}}
 \put(30,70){\line(1,0){40}} \put(70,40){\line(0,1){30}}

\put(70,70){\line(-5,1){50}} \put(17,76){$\blacksquare$}
\put(0,10){\line(1,4){18}}

  \put(0,10){\line(1,6){12}} \put(100,100){\line(-6,-1){90}}
 \put(7,80){$\blacksquare$}

 \put(20,-18){Figure 2(iii)(v)}

\end{picture}
\end{center}
$$$$\\\\\\\\\\\\

\begin{center}
\begin{picture}(0,0)(165,104)

\put(-2,8){$\bullet$} \put(-2,-2){$a_1$}

\put(-2,98){$\blacksquare$} \put(-2,108){$b_1$}

\put(98,8){$\blacksquare$} \put(98,-2){$b_2$}

\put(98,98){$\bullet$} \put(98,104){$a_2$}

\put(0,10){\line(1,0){100}} \put(0,10){\line(0,1){90}}
\put(0,100){\line(1,0){100}} \put(100,10){\line(0,1){90}}

\put(28,38){$\bullet$} \put(33,43){$a_1$}

\put(68,38){$\blacksquare$} \put(60,43){$b_2$}

\put(28,68){$\blacksquare$} \put(37,62){$b_1$}

\put(68,68){$\bullet$} \put(60,63){$a_2$}

 \put(30,40){\line(1,0){40}}  \put(30,40){\line(0,1){30}}
 \put(30,70){\line(1,0){40}} \put(70,40){\line(0,1){30}}

  \put(0,10){\line(1,6){13}} \put(70,70){\line(-3,1){60}}
 \put(8,86){$\blacksquare$}

\put(30,40){\line(-1,4){10}} \put(70,70){\line(-4,1){50}}
\put(16,77){$\blacksquare$}

 \put(20,-18){Figure 2(iii)(vi)}

\end{picture}
\begin{picture}(0,0)(50,104)

\put(-2,8){$\bullet$} \put(-2,-2){$a_1$}

\put(-2,98){$\blacksquare$} \put(-2,108){$b_1$}

\put(98,8){$\blacksquare$} \put(98,-2){$b_2$}

\put(98,98){$\bullet$} \put(98,104){$a_2$}

\put(0,10){\line(1,0){100}} \put(0,10){\line(0,1){90}}
\put(0,100){\line(1,0){100}} \put(100,10){\line(0,1){90}}

\put(28,38){$\bullet$} \put(33,43){$a_1$}

\put(68,38){$\blacksquare$} \put(60,43){$b_2$}

\put(28,68){$\blacksquare$} \put(37,62){$b_1$}

\put(68,68){$\bullet$} \put(60,63){$a_2$}

 \put(30,40){\line(1,0){40}}  \put(30,40){\line(0,1){30}}
 \put(30,70){\line(1,0){40}} \put(70,40){\line(0,1){30}}

\put(30,40){\line(4,-1){50}} \put(71,70){\line(1,-4){10}}
\put(77,27){$\blacksquare$}

  \put(0,10){\line(1,6){13}} \put(70,70){\line(-3,1){60}}
 \put(8,86){$\blacksquare$}

 \put(20,-18){Figure 2(iii)(vii)}

\end{picture}
\begin{picture}(0,0)(-65,104)

\put(-2,8){$\bullet$} \put(-2,-2){$a_1$}

\put(-2,98){$\blacksquare$} \put(-2,108){$b_1$}

\put(98,8){$\blacksquare$} \put(98,-2){$b_2$}

\put(98,98){$\bullet$} \put(98,104){$a_2$}

\put(0,10){\line(1,0){100}} \put(0,10){\line(0,1){90}}
\put(0,100){\line(1,0){100}} \put(100,10){\line(0,1){90}}

\put(28,38){$\bullet$} \put(33,43){$a_1$}

\put(68,38){$\blacksquare$} \put(60,43){$b_2$}

\put(28,68){$\blacksquare$} \put(37,62){$b_1$}

\put(68,68){$\bullet$} \put(60,63){$a_2$}

 \put(30,40){\line(1,0){40}}  \put(30,40){\line(0,1){30}}
 \put(30,70){\line(1,0){40}} \put(70,40){\line(0,1){30}}

  \put(0,10){\line(1,6){13}} \put(70,70){\line(-3,1){60}}
 \put(8,86){$\blacksquare$}

 \put(100,100){\line(-1,-6){12}} \put(85,25){$\blacksquare$}
 \put(87,30){\line(-6,1){57}}

 \put(20,-18){Figure 2(iii)(viii)}

\end{picture}
\end{center}
$$$$\\\\\\\\\\\\

\begin{center}
\begin{picture}(0,0)(165,104)

\put(-2,8){$\bullet$} \put(-2,-2){$a_1$}

\put(-2,98){$\blacksquare$} \put(-2,108){$b_1$}

\put(98,8){$\blacksquare$} \put(98,-2){$b_2$}

\put(98,98){$\bullet$} \put(98,104){$a_2$}

\put(0,10){\line(1,0){100}} \put(0,10){\line(0,1){90}}
\put(0,100){\line(1,0){100}} \put(100,10){\line(0,1){90}}

\put(28,38){$\bullet$} \put(33,43){$a_1$}

\put(68,38){$\blacksquare$} \put(60,43){$b_2$}

\put(28,68){$\blacksquare$} \put(37,62){$b_1$}

\put(68,68){$\bullet$} \put(60,63){$a_2$}

 \put(30,40){\line(1,0){40}}  \put(30,40){\line(0,1){30}}
 \put(30,70){\line(1,0){40}} \put(70,40){\line(0,1){30}}

 \put(0,10){\line(6,1){90}}
 \put(100,100){\line(-1,-6){12}} \put(85,25){$\blacksquare$}

 \put(0,10){\line(3,1){85}}
 \put(100,100){\line(-1,-4){15}} \put(79,35){$\blacksquare$}

 \put(10,-18){Figure 2(iv)(iv)}

\end{picture}
\begin{picture}(0,0)(50,104)

\put(-2,8){$\bullet$} \put(-2,-2){$a_1$}

\put(-2,98){$\blacksquare$} \put(-2,108){$b_1$}

\put(98,8){$\blacksquare$} \put(98,-2){$b_2$}

\put(98,98){$\bullet$} \put(98,104){$a_2$}

\put(0,10){\line(1,0){100}} \put(0,10){\line(0,1){90}}
\put(0,100){\line(1,0){100}} \put(100,10){\line(0,1){90}}

\put(28,38){$\bullet$} \put(33,43){$a_1$}

\put(68,38){$\blacksquare$} \put(60,43){$b_2$}

\put(28,68){$\blacksquare$} \put(37,62){$b_1$}

\put(68,68){$\bullet$} \put(60,63){$a_2$}

 \put(30,40){\line(1,0){40}}  \put(30,40){\line(0,1){30}}
 \put(30,70){\line(1,0){40}} \put(70,40){\line(0,1){30}}

 \put(0,10){\line(6,1){90}}
 \put(100,100){\line(-1,-6){12}} \put(85,25){$\blacksquare$}

 \put(0,10){\line(1,6){12}} \put(100,100){\line(-6,-1){90}}
 \put(7,80){$\blacksquare$}

 \put(10,-18){Figure 2(iv)(v)}

\end{picture}
\begin{picture}(0,0)(-65,104)

\put(-2,8){$\bullet$} \put(-2,-2){$a_1$}

\put(-2,98){$\blacksquare$} \put(-2,108){$b_1$}

\put(98,8){$\blacksquare$} \put(98,-2){$b_2$}

\put(98,98){$\bullet$} \put(98,104){$a_2$}

\put(0,10){\line(1,0){100}} \put(0,10){\line(0,1){90}}
\put(0,100){\line(1,0){100}} \put(100,10){\line(0,1){90}}

\put(28,38){$\bullet$} \put(33,43){$a_1$}

\put(68,38){$\blacksquare$} \put(60,43){$b_2$}

\put(28,68){$\blacksquare$} \put(37,62){$b_1$}

\put(68,68){$\bullet$} \put(60,63){$a_2$}

 \put(30,40){\line(1,0){40}}  \put(30,40){\line(0,1){30}}
 \put(30,70){\line(1,0){40}} \put(70,40){\line(0,1){30}}

 \put(0,10){\line(6,1){90}}
 \put(100,100){\line(-1,-6){12}} \put(85,25){$\blacksquare$}

 \put(30,40){\line(-1,4){10}} \put(70,70){\line(-4,1){50}}
\put(16,78){$\blacksquare$}

 \put(10,-18){Figure 2(iv)(vi)}

\end{picture}
\end{center}
$$$$\\\\\\\\\\\\

\newpage

\begin{center}
\begin{picture}(0,0)(120,104)

\put(-2,8){$\bullet$} \put(-2,-2){$a_1$}

\put(-2,98){$\blacksquare$} \put(-2,108){$b_1$}

\put(98,8){$\blacksquare$} \put(98,-2){$b_2$}

\put(98,98){$\bullet$} \put(98,104){$a_2$}

\put(0,10){\line(1,0){100}} \put(0,10){\line(0,1){90}}
\put(0,100){\line(1,0){100}} \put(100,10){\line(0,1){90}}

\put(28,38){$\bullet$} \put(33,43){$a_1$}

\put(68,38){$\blacksquare$} \put(60,43){$b_2$}

\put(28,68){$\blacksquare$} \put(37,62){$b_1$}

\put(68,68){$\bullet$} \put(60,63){$a_2$}

 \put(30,40){\line(1,0){40}}  \put(30,40){\line(0,1){30}}
 \put(30,70){\line(1,0){40}} \put(70,40){\line(0,1){30}}

 \put(0,10){\line(6,1){90}}
 \put(100,100){\line(-1,-6){12}} \put(85,25){$\blacksquare$}

\put(30,40){\line(4,-1){50}} \put(71,70){\line(1,-4){10}}
\put(77,27){$\blacksquare$}

 \put(10,-18){Figure 2(iv)(vii)}

\end{picture}
\begin{picture}(0,0)(-25,104)

\put(-2,8){$\bullet$} \put(-2,-2){$a_1$}

\put(-2,98){$\blacksquare$} \put(-2,108){$b_1$}

\put(98,8){$\blacksquare$} \put(98,-2){$b_2$}

\put(98,98){$\bullet$} \put(98,104){$a_2$}

\put(0,10){\line(1,0){100}} \put(0,10){\line(0,1){90}}
\put(0,100){\line(1,0){100}} \put(100,10){\line(0,1){90}}

\put(28,38){$\bullet$} \put(33,43){$a_1$}

\put(68,38){$\blacksquare$} \put(60,43){$b_2$}

\put(28,68){$\blacksquare$} \put(37,62){$b_1$}

\put(68,68){$\bullet$} \put(60,63){$a_2$}

 \put(30,40){\line(1,0){40}}  \put(30,40){\line(0,1){30}}
 \put(30,70){\line(1,0){40}} \put(70,40){\line(0,1){30}}

 \put(0,10){\line(6,1){90}}
 \put(100,100){\line(-1,-6){12}} \put(85,25){$\blacksquare$}

 \put(100,100){\line(-1,-3){23}}
 \put(77,32){\line(-6,1){48}} \put(75,29){$\blacksquare$}

 \put(10,-18){Figure 2(iv)(viii)}

\end{picture}
\end{center}
$$$$\\\\\\\\\\\\

\begin{center}
\begin{picture}(0,0)(120,104)

\put(-2,8){$\bullet$} \put(-2,-2){$a_1$}

\put(-2,98){$\blacksquare$} \put(-2,108){$b_1$}

\put(98,8){$\blacksquare$} \put(98,-2){$b_2$}

\put(98,98){$\bullet$} \put(98,104){$a_2$}

\put(0,10){\line(1,0){100}} \put(0,10){\line(0,1){90}}
\put(0,100){\line(1,0){100}} \put(100,10){\line(0,1){90}}

\put(28,38){$\bullet$} \put(33,43){$a_1$}

\put(68,38){$\blacksquare$} \put(60,43){$b_2$}

\put(28,68){$\blacksquare$} \put(37,62){$b_1$}

\put(68,68){$\bullet$} \put(60,63){$a_2$}

 \put(30,40){\line(1,0){40}}  \put(30,40){\line(0,1){30}}
 \put(30,70){\line(1,0){40}} \put(70,40){\line(0,1){30}}

  \put(0,10){\line(1,6){12}} \put(100,100){\line(-6,-1){90}}
 \put(7,80){$\blacksquare$}

 \put(100,100){\line(-4,-1){80}} \put(0,10){\line(1,3){23}}
 \put(17,75){$\blacksquare$}

 \put(10,-18){Figure 2(v)(v)}

\end{picture}
\begin{picture}(0,0)(-25,104)

\put(-2,8){$\bullet$} \put(-2,-2){$a_1$}

\put(-2,98){$\blacksquare$} \put(-2,108){$b_1$}

\put(98,8){$\blacksquare$} \put(98,-2){$b_2$}

\put(98,98){$\bullet$} \put(98,104){$a_2$}

\put(0,10){\line(1,0){100}} \put(0,10){\line(0,1){90}}
\put(0,100){\line(1,0){100}} \put(100,10){\line(0,1){90}}

\put(28,38){$\bullet$} \put(33,43){$a_1$}

\put(68,38){$\blacksquare$} \put(60,43){$b_2$}

\put(28,68){$\blacksquare$} \put(37,62){$b_1$}

\put(68,68){$\bullet$} \put(60,63){$a_2$}

 \put(30,40){\line(1,0){40}}  \put(30,40){\line(0,1){30}}
 \put(30,70){\line(1,0){40}} \put(70,40){\line(0,1){30}}

 \put(0,10){\line(1,6){12}} \put(100,100){\line(-6,-1){90}}
 \put(7,80){$\blacksquare$}

 \put(30,40){\line(-1,4){10}} \put(70,70){\line(-4,1){50}}
\put(16,78){$\blacksquare$}

 \put(10,-18){Figure 2(v)(vi)}

\end{picture}
\end{center}
$$$$\\\\\\\\\\\\

\begin{center}
\begin{picture}(0,0)(120,104)

\put(-2,8){$\bullet$} \put(-2,-2){$a_1$}

\put(-2,98){$\blacksquare$} \put(-2,108){$b_1$}

\put(98,8){$\blacksquare$} \put(98,-2){$b_2$}

\put(98,98){$\bullet$} \put(98,104){$a_2$}

\put(0,10){\line(1,0){100}} \put(0,10){\line(0,1){90}}
\put(0,100){\line(1,0){100}} \put(100,10){\line(0,1){90}}

\put(28,38){$\bullet$} \put(33,43){$a_1$}

\put(68,38){$\blacksquare$} \put(60,43){$b_2$}

\put(28,68){$\blacksquare$} \put(37,62){$b_1$}

\put(68,68){$\bullet$} \put(60,63){$a_2$}

 \put(30,40){\line(1,0){40}}  \put(30,40){\line(0,1){30}}
 \put(30,70){\line(1,0){40}} \put(70,40){\line(0,1){30}}

 \put(0,10){\line(1,6){12}} \put(100,100){\line(-6,-1){90}}
 \put(7,80){$\blacksquare$}

\put(30,40){\line(4,-1){50}} \put(71,70){\line(1,-4){10}}
\put(77,27){$\blacksquare$}

 \put(10,-18){Figure 2(v)(vii)}

\end{picture}
\begin{picture}(0,0)(-25,104)

\put(-2,8){$\bullet$} \put(-2,-2){$a_1$}

\put(-2,98){$\blacksquare$} \put(-2,108){$b_1$}

\put(98,8){$\blacksquare$} \put(98,-2){$b_2$}

\put(98,98){$\bullet$} \put(98,104){$a_2$}

\put(0,10){\line(1,0){100}} \put(0,10){\line(0,1){90}}
\put(0,100){\line(1,0){100}} \put(100,10){\line(0,1){90}}

\put(28,38){$\bullet$} \put(33,43){$a_1$}

\put(68,38){$\blacksquare$} \put(60,43){$b_2$}

\put(28,68){$\blacksquare$} \put(37,62){$b_1$}

\put(68,68){$\bullet$} \put(60,63){$a_2$}

 \put(30,40){\line(1,0){40}}  \put(30,40){\line(0,1){30}}
 \put(30,70){\line(1,0){40}} \put(70,40){\line(0,1){30}}

 \put(0,10){\line(1,6){12}} \put(100,100){\line(-6,-1){90}}
 \put(7,80){$\blacksquare$}

 \put(100,100){\line(-1,-3){23}}
 \put(77,32){\line(-6,1){48}} \put(75,29){$\blacksquare$}

 \put(10,-18){Figure 2(v)(viii)}

\end{picture}
\end{center}
$$$$\\\\\\\\\\\\

\newpage

\begin{center}
\begin{picture}(0,0)(165,104)

\put(-2,8){$\bullet$} \put(-2,-2){$a_1$}

\put(-2,98){$\blacksquare$} \put(-2,108){$b_1$}

\put(98,8){$\blacksquare$} \put(98,-2){$b_2$}

\put(98,98){$\bullet$} \put(98,104){$a_2$}

\put(0,10){\line(1,0){100}} \put(0,10){\line(0,1){90}}
\put(0,100){\line(1,0){100}} \put(100,10){\line(0,1){90}}

\put(28,38){$\bullet$} \put(33,43){$a_1$}

\put(68,38){$\blacksquare$} \put(60,43){$b_2$}

\put(28,68){$\blacksquare$} \put(37,62){$b_1$}

\put(68,68){$\bullet$} \put(60,63){$a_2$}

 \put(30,40){\line(1,0){40}}  \put(30,40){\line(0,1){30}}
 \put(30,70){\line(1,0){40}} \put(70,40){\line(0,1){30}}

\put(30,40){\line(-1,4){10}} \put(70,70){\line(-4,1){50}}
\put(16,78){$\blacksquare$}

\put(30,40){\line(-1,2){25}} \put(70,70){\line(-3,2){30}}
\put(3,86){$\blacksquare$} \put(40,90){\line(-1,0){30}}

 \put(10,-18){Figure 2(vi)(vi)}

\end{picture}
\begin{picture}(0,0)(50,104)

\put(-2,8){$\bullet$} \put(-2,-2){$a_1$}

\put(-2,98){$\blacksquare$} \put(-2,108){$b_1$}

\put(98,8){$\blacksquare$} \put(98,-2){$b_2$}

\put(98,98){$\bullet$} \put(98,104){$a_2$}

\put(0,10){\line(1,0){100}} \put(0,10){\line(0,1){90}}
\put(0,100){\line(1,0){100}} \put(100,10){\line(0,1){90}}

\put(28,38){$\bullet$} \put(33,43){$a_1$}

\put(68,38){$\blacksquare$} \put(60,43){$b_2$}

\put(28,68){$\blacksquare$} \put(37,62){$b_1$}

\put(68,68){$\bullet$} \put(60,63){$a_2$}

 \put(30,40){\line(1,0){40}}  \put(30,40){\line(0,1){30}}
 \put(30,70){\line(1,0){40}} \put(70,40){\line(0,1){30}}

\put(30,40){\line(-1,4){10}} \put(70,70){\line(-4,1){50}}
\put(16,78){$\blacksquare$}

\put(30,40){\line(4,-1){50}} \put(71,70){\line(1,-4){10}}
\put(77,27){$\blacksquare$}

 \put(10,-18){Figure 2(vi)(vii)}

\end{picture}
\begin{picture}(0,0)(-65,104)

\put(-2,8){$\bullet$} \put(-2,-2){$a_1$}

\put(-2,98){$\blacksquare$} \put(-2,108){$b_1$}

\put(98,8){$\blacksquare$} \put(98,-2){$b_2$}

\put(98,98){$\bullet$} \put(98,104){$a_2$}

\put(0,10){\line(1,0){100}} \put(0,10){\line(0,1){90}}
\put(0,100){\line(1,0){100}} \put(100,10){\line(0,1){90}}

\put(28,38){$\bullet$} \put(33,43){$a_1$}

\put(68,38){$\blacksquare$} \put(60,43){$b_2$}

\put(28,68){$\blacksquare$} \put(37,62){$b_1$}

\put(68,68){$\bullet$} \put(60,63){$a_2$}

 \put(30,40){\line(1,0){40}}  \put(30,40){\line(0,1){30}}
 \put(30,70){\line(1,0){40}} \put(70,40){\line(0,1){30}}

\put(30,40){\line(-1,4){10}} \put(70,70){\line(-4,1){50}}
\put(16,78){$\blacksquare$}

 \put(100,100){\line(-1,-6){12}} \put(85,25){$\blacksquare$}
 \put(87,30){\line(-6,1){57}}

 \put(10,-18){Figure 2(vi)(viii)}

\end{picture}
\end{center}
$$$$\\\\\\\\\\\\

\begin{center}
\begin{picture}(0,0)(120,104)

\put(-2,8){$\bullet$} \put(-2,-2){$a_1$}

\put(-2,98){$\blacksquare$} \put(-2,108){$b_1$}

\put(98,8){$\blacksquare$} \put(98,-2){$b_2$}

\put(98,98){$\bullet$} \put(98,104){$a_2$}

\put(0,10){\line(1,0){100}} \put(0,10){\line(0,1){90}}
\put(0,100){\line(1,0){100}} \put(100,10){\line(0,1){90}}

\put(28,38){$\bullet$} \put(33,43){$a_1$}

\put(68,38){$\blacksquare$} \put(60,43){$b_2$}

\put(28,68){$\blacksquare$} \put(37,62){$b_1$}

\put(68,68){$\bullet$} \put(60,63){$a_2$}

 \put(30,40){\line(1,0){40}}  \put(30,40){\line(0,1){30}}
 \put(30,70){\line(1,0){40}} \put(70,40){\line(0,1){30}}

\put(30,40){\line(4,-1){50}} \put(71,70){\line(1,-4){10}}
\put(77,27){$\blacksquare$}

\put(30,40){\line(3,-1){60}} \put(71,70){\line(1,-2){23}}
\put(88,18){$\blacksquare$}

 \put(10,-18){Figure 2(vii)(vii)}

\end{picture}
\begin{picture}(0,0)(-25,104)

\put(-2,8){$\bullet$} \put(-2,-2){$a_1$}

\put(-2,98){$\blacksquare$} \put(-2,108){$b_1$}

\put(98,8){$\blacksquare$} \put(98,-2){$b_2$}

\put(98,98){$\bullet$} \put(98,104){$a_2$}

\put(0,10){\line(1,0){100}} \put(0,10){\line(0,1){90}}
\put(0,100){\line(1,0){100}} \put(100,10){\line(0,1){90}}

\put(28,38){$\bullet$} \put(33,43){$a_1$}

\put(68,38){$\blacksquare$} \put(60,43){$b_2$}

\put(28,68){$\blacksquare$} \put(37,62){$b_1$}

\put(68,68){$\bullet$} \put(60,63){$a_2$}

 \put(30,40){\line(1,0){40}}  \put(30,40){\line(0,1){30}}
 \put(30,70){\line(1,0){40}} \put(70,40){\line(0,1){30}}

\put(30,40){\line(4,-1){50}} \put(71,70){\line(1,-4){10}}
\put(77,27){$\blacksquare$}

\put(30,40){\line(3,-1){58}} \put(99,100){\line(-1,-6){13}}
\put(83,18){$\blacksquare$}

 \put(10,-18){Figure 2(vii)(viii)}

\end{picture}
\end{center}
$$$$\\\\\\\\\\\\

\begin{center}
\begin{picture}(0,0)(50,104)

\put(-2,8){$\bullet$} \put(-2,-2){$a_1$}

\put(-2,98){$\blacksquare$} \put(-2,108){$b_1$}

\put(98,8){$\blacksquare$} \put(98,-2){$b_2$}

\put(98,98){$\bullet$} \put(98,104){$a_2$}

\put(0,10){\line(1,0){100}} \put(0,10){\line(0,1){90}}
\put(0,100){\line(1,0){100}} \put(100,10){\line(0,1){90}}

\put(28,38){$\bullet$} \put(33,43){$a_1$}

\put(68,38){$\blacksquare$} \put(60,43){$b_2$}

\put(28,68){$\blacksquare$} \put(37,62){$b_1$}

\put(68,68){$\bullet$} \put(60,63){$a_2$}

 \put(30,40){\line(1,0){40}}  \put(30,40){\line(0,1){30}}
 \put(30,70){\line(1,0){40}} \put(70,40){\line(0,1){30}}

\put(30,40){\line(4,-1){50}} \put(99,100){\line(-1,-4){17}}
\put(77,27){$\blacksquare$}

\put(30,40){\line(3,-1){58}} \put(99,100){\line(-1,-6){13}}
\put(83,18){$\blacksquare$}

 \put(10,-18){Figure 2(viii)(viii)}

\end{picture}
\end{center}
$$$$\\\\\\\\\\\\\\

The following tables contain the following information: in each
table, the Figure in the first row is isomorphic to the Figure in
the second row which appeared in \cite{Tung}. If ``Not included"
appears in the second row, it means that the corresponding Figure in
the first row does not have a region which contains all $b_i$,
therefore it is not included. For example,  Figure 2(i)(i) is
isomorphic to Figure 7(ii) in \cite{Tung}; Figure 2(i)(ii) is
isomorphic to Figure 7(i) in \cite{Tung}; Figure 2(i)(vi) is
isomorphic to a drawing symmetric to Figure 7(i) in \cite{Tung};
Figure 2(iv)(iv) is not included since it does not have a region
which contains all $b_i$. These tables show that Figures 7(i) to
7(vii) in \cite{Tung} exhaust
all the possible drawings.\\\\
\begin{tabular}{|c|c|c|c|c|c|c|}
  \hline
 2(i)(i) &  2(i)(ii) &  2(i)(iv) & 2(i)(v) &  2(i)(vi) &  2(i)(vii)  &  2(i)(viii)\\
  \hline
  7(ii)  & 7(i) & 7(iii) & sym 7(iii)  & sym 7(vi) &  7(vi) & 7(viii) \\
  \hline
\end{tabular}
\\\\\\
\begin{tabular}{|c|c|c|c|c|c|}
  \hline
 2(ii)(ii) &  2(ii)(iii) &  2(ii)(iv) & 2(ii)(v) &  2(ii)(vi) &  2(ii)(vii)  \\
  \hline
  7(ii)  & 7(viii) & sym 7(iii) &  7(iii)  & 7(vi) & sym 7(vi) \\
  \hline
\end{tabular}
\\\\\\
\begin{tabular}{|c|c|c|c|c|c|}
  \hline
2(iii)(iii) &  2(iii)(iv) & 2(iii)(v) &  2(iii)(vi) &  2(iii)(vii) &  2(iii)(viii)  \\
  \hline
  7(ii)  & sym 7(iii) & 7(iii) & 7(vi)  & sym 7(vi) & sym 7(i) \\
  \hline
\end{tabular}
\\\\\\
\begin{tabular}{|c|c|c|c|c|c|}
  \hline
  2(iv)(iv) & 2(iv)(v) &  2(iv)(vi) &  2(iv)(vii) &  2(iv)(viii)  \\
  \hline
  Not included  & 7(v) & 7(iv) & Not included  & 7(iii)  \\
  \hline
\end{tabular}
\\\\\\
\begin{tabular}{|c|c|c|c|c|c|}
  \hline
 2(v)(v) &  2(v)(vi) &  2(v)(vii) &  2(v)(viii)  \\
  \hline
  Not included & Not included & 7(iv)  & sym 7(iii)  \\
  \hline
\end{tabular}
\\\\\\
\begin{tabular}{|c|c|c|c|c|c|}
  \hline
 2(vi)(vi) &  2(vi)(vii) &  2(vi)(viii)  \\
  \hline
   Not included & 7(vii)  & sym 7(vi)  \\
  \hline
\end{tabular}
\\\\\\
\begin{tabular}{|c|c|c|c|c|c|}
  \hline
2(vii)(vii) &  2(vii)(viii)  \\
  \hline
   Not included &  7(vi)  \\
  \hline
\end{tabular}
\\\\\\
\begin{tabular}{|c|c|c|c|c|c|}
  \hline
2(viii)(viii)  \\
  \hline
 7(ii)  \\
  \hline
\end{tabular}

\end{document}